\documentclass[11pt]{article}
\usepackage{amsmath,amssymb,amsthm,color}

\usepackage{comment}
\usepackage{color}
\usepackage{amssymb}
\usepackage{amsmath}
\usepackage[english]{babel}
\usepackage{fancyhdr}
\usepackage{amsmath}

\newtheorem{thm}{Theorem}[section]

\newtheorem{lemma}{Lemma}[section]

\newtheorem{cor}{Corollary}[section]

\def\ep{\hfill $\Diamond$}
\def\R{{\mathfrak R}\, }

\def\ci{\begin{color}{red}\,}
\def\cf{\end{color}\,}

\begin{document}
\begin{center}
{\Large\bf Generalized Jordan derivations of Incidence Algebras}

\vspace{.2in}
{\bf Bruno Leonardo Macedo Ferreira\textsuperscript{1}} 
\\
{\bf Tanise Carnieri Pierin\textsuperscript{2}}
\\
and
\\
{\bf Ruth Nascimento Ferreira\textsuperscript{3}}
\vspace{.2in}

 Universidade Tecnol\'{o}gica Federal do Paran\'{a}\textsuperscript{1,3}, Av. Profa. Laura Pacheco Bastos, 800, 85053-510, Guarapuava, Brazil.
\\
and
\\
 Universidade Federal do Paran\'{a}\textsuperscript{2}, Av. Cel. Francisco H. dos Santos, 100, 81530-000, Curitiba, Brazil.
\vspace{.2in}
 
brunoferreira@utfpr.edu.br
\\
tanise@ufpr.br
\\
and
\\
ruthnascimento@utfpr.edu.br
\vspace{.2in}

\end{center}

\begin{abstract}
For a given ring $\R$ and a locally finite pre-ordered set $(X, \leq)$, consider $I(X, \R)$ to be the incidence algebra of $X$ over $\R$. 
Motivated by a Xiao's result which states that every Jordan derivation of $I(X,\R)$ is a derivation in the case $\R$ is $2$-torsion free, one proves that each generalized Jordan derivation of $I(X,\R)$ is a generalized derivation provided $\R$ is $2$-torsion free, getting as a consequence the above mentioned result.
\end{abstract}

\noindent {\bf Mathematics Subject Classification (2010):} 16W25; 47B47.

\noindent {\bf Keywords:} generalized derivation; generalized Jordan derivation; incidence algebras.

\section{Introduction}
For a given ring $\R$, recall that a linear map $d$ from $\R$ into itself is called a derivation if $d(ab) = d(a)b+ad(b)$
for all $a,b \in \R$; and a Jordan derivation if $d(a^2) = d(a)a+ad(a)$ for each $a \in \R$. More generally \cite{Jing}, if
there is a derivation $\tau : \R \rightarrow \R$ such that $d(ab) = d(a)b +a\tau (b)$ for all $a,b \in \R$, then $d$ is
called a generalized derivation and $\tau$ is the relating derivation; analogously, if there is a Jordan derivation
$\tau : \R \rightarrow \R$ such that $d(a^2) = d(a)a + a\tau(a)$ for all $a \in \R$, then $d$ is called a generalized
Jordan derivation and $\tau$ is the relating Jordan derivation. 
The structures of derivations, Jordan derivations, generalized derivations and generalized
Jordan derivations were systematically studied. It is obvious that every generalized derivation
is a generalized Jordan derivation and every derivation is a Jordan derivation. But the converse is in general not true.
Herstein \cite{Herstein} showed that every Jordan derivation from a $2$-torsion free prime ring into itself is a derivation.
Bre$\check{s}$ar \cite{Bresar} proved that Herstein's result is true for $2$-torsion free semiprime rings.
Jing and Lu, motivated by the concept of generalized derivation, introduce this concept of generalized Jordan derivation in \cite{Jing}. 

Let us now recall the notion of incidence algebra \cite{kopp}, \cite{spi}, which we deal
in this paper. Let $(X,\leq)$ be a locally finite
pre-ordered set. This means $\leq$ is a reflexive and transitive binary relation on the
set $X$, and for any $x \leq y$ in $X$ there are only finitely many elements $z$ satisfying
$x \leq z \leq y$. The incidence algebra $I(X,\R)$ of $X$ over $\R$ is defined as the set
$$I(X,\R):= \left\{f : X \times X \rightarrow \R ~ | ~ f(x, y) = 0 ~ if ~ x \nleq y\right\}$$
with algebraic operation given by
$$(f + g)(x, y) = f(x, y) + g(x, y),$$
$$(rf)(x, y) = rf(x, y),$$
$$(fg)(x, y) = \sum_{x\leq z \leq y} f(x, z)g(z, y)$$
for all $f, g \in I(X,\R)$, $r \in \R$ and $x, y, z \in X$. The product $fg$ is usually called
convolution in function theory. It would be helpful to point out that the full matrix
algebra $M_n(\R)$ and the upper (or lower) triangular matrix algebras $T_n(\R)$ are
special examples of incidence algebras.
The identity element $\delta$ of $I(X,\R)$ is given by $\delta(x,y) = \delta_{xy}$ for $x \leq y$, where $\delta_{xy} \in \{0, 1\}$ is the Kronecker delta. For given $x, y \in X$ with $x \leq y$, let $e_{xy}$ be defined by $e_{xy}(u, v) = 1$ if $(u, v) = (x, y)$, and $e_{xy}(u, v) = 0$
otherwise. Then $e_{xy}e_{uv} = \delta_{yu}e_{xv}$ by the definition of convolution. Moreover, the set $B := \left\{e_{xy} | x \leq y\right\}$ forms an $\R$-linear basis of
$I(X,\R)$. Note that incidence algebras allow infinite summation, and hence the $\R$-linear map here means a map preserving infinite sum and scalar multiplication.

Incidence algebras were first considered by Ward \cite{ward} as generalized algebras of arithmetic functions. Rota and Stanley developed incidence algebras as the fundamental structures of enumerative combinatorial theory and allied areas of arithmetic function theory (see \cite{spiegel}). Motivated by the results of Stanley \cite{stanley}, automorphisms and other algebraic mappings of incidence algebras have been extensively studied (see \cite{bacla}, \cite{coelho}, \cite{jøndrup}, \cite{kopp}, \cite{mathis}, \cite{nowicki}, \cite{nowicki2}, \cite{spiegel} and the references therein). Baclawski \cite{bacla} studied the automorphisms and derivations of incidence algebras $I(X,\R)$ when $X$ is a locally finite partially ordered set. More specifically, he proved that every derivation of $I(X,\R)$ with $X$ a locally finite partially ordered set can be decomposed as a sum of an inner derivation and a transitive induced derivation. Koppinen \cite{kopp} has extended these results to the incidence algebras $I(X,\R)$ with $X$ a locally finite pre-ordered set. Xiao \cite{Xiao} proved that every Jordan derivation of $I(X,\R)$ is a derivation provided that $\R$ is $2$-torsion free.
Motivated by Xiao's result our main objective is to prove that every generalized Jordan derivation of $I(X,\R)$ is a generalized derivation provided that $\R$ is $2$-torsion free.
 
\section{Results}

We first collect some background material to prove our main result. 
Throughout this section, $\R$ denotes a  $2$-torsion free ring. Let $\Xi: I(X,\R) \rightarrow I(X,\R)$ be a
generalized Jordan derivation and $\tau : I(X,\R) \rightarrow I(X,\R)$ the relating Jordan derivation.

\begin{lemma}\label{lema1} For all $a, b, c \in I(X,\R)$, the following statements hold:
\begin{enumerate}
\item[(1)] $\Xi(ab + ba) = \Xi(a)b + a\tau (b) + \Xi(b)a + b\tau (a),$
\item[(2)] $\Xi(aba) = \Xi(a)ba + a\tau (b)a + ab\tau (a),$
\item[(3)] $\Xi(abc + cba) = \Xi(a)bc + a\tau (b)c + ab\tau (c) + \Xi(c)ba + c\tau (b)a + cb\tau (a).$
\end{enumerate}
\end{lemma}
\proof See \cite{Jing}. \ep

\vspace{0.5cm}

According to Lemma \ref{lema1}, $\Xi(aba) = \Xi(a)ba + a\tau (b)a + ab\tau (a)$. In the case $ab =ba =0$, we obtain $a\tau(b)a =0$. Furthermore, it follows that
\begin{eqnarray}\label{id1}
\Xi(e) =\Xi(e)e + e\tau(e),
\end{eqnarray}
for any idempotent $e \in I(X,\R)$. In particular, since (\ref{id1}), $e\tau(a)e =0$, for any $a \in I(X,\R)$
satisfying $ea =ae =0$, and $\Xi(a)e + a\tau(e) + \Xi(e)a + e\tau(a)=0$. Multiplying by $e$ on the right yields
\begin{eqnarray}\label{id2}
\Xi(a)e + a\tau(e)=0 =\Xi(e)a + e\tau(a),
\end{eqnarray}
for any idempotent $e$ satisfying $ea =ae =0$.

Now assume that the set $B := \left\{e_{xy} | x \leq y\right\}$ forms an $\R$-linear basis of $I(X,\R)$. It is a consequence of (\ref{id1}) that
\begin{eqnarray}\label{id3}
\Xi(e_{ii}) = \Xi(e_{ii})e_{ii} + e_{ii}\tau(e_{ii})
 \ \ \ \ \ \ \mbox{and} \ \ \ \ \ \
e_{ki}\tau(e_{ii})e_{ij} = 0,
\end{eqnarray}
for all $i$ and $k \leq i \leq j$. From Lemma \ref{lema1} and the fact that $\Xi(e_{ij}) = \Xi(e_{ii}e_{ij} + e_{ij}e_{ii})$ for all
$1 \leq i < j \leq n$, we obtain
\begin{eqnarray}\label{id4}
\Xi(e_{ij}) = \Xi(e_{ii})e_{ij} + e_{ii}\tau(e_{ij}) + \Xi(e_{ij})e_{ii} + e_{ij}\tau(e_{ii}) 
\end{eqnarray}
whenever $i < j$. Furthermore (\ref{id2}) implies that
\begin{eqnarray}\label{id5}
\Xi(e_{kj})e_{ii} + e_{kj}\Xi(e_{ii}) = \Xi(e_{ii})e_{kj} + e_{ii}\tau(e_{kj}) = 0 
\end{eqnarray}
for all $k, j \neq i$.
Define a $\R$-linear map $\phi$ from $I(X,\R)$ into itself by letting
\begin{eqnarray}\label{id6}
\phi(e_{ij}) = \Xi(e_{ii})e_{ij} + e_{ii}\tau(e_{ij}), \ \ \ \ \ \ i\leq j. 
\end{eqnarray}
According to (\ref{id3}), $\phi(e_{ii}) = \Xi(e_{ii})$.
Xiao proved the following result. 
\begin{lemma}[Lemma 3.2 \cite{Xiao}]\label{xiao1} 
Let $\tau : I(X,\R) \rightarrow I(X,\R)$ be a Jordan derivation. Then
$$\tau(e_{ij}) = \sum_{x \in L_i}C^{ii}_{xi}e_{xj} + C_{ij}^{ij}e_{ij} + \sum_{y \in R_j}C^{jj}_{jy}e_{iy} + C^{ij}_{ji} e_{ji}$$
for all $e_{ij} \in B$, where the coefficients $C^{ij}_{xy}$ are subject to the following relations
$$C^{jj}_{jk} + C^{kk}_{jk} = 0,  \ \ \ \ if j \leq k;$$
$$C^{ij}_{ij} + C_{jk}^{jk} = C_{ik}^{ik} , \ \ \ \ \ if i \leq j, j \leq k.$$
\end{lemma} 

\begin{lemma}\label{lema2}
$\phi$ is a generalized derivation.
\end{lemma}
\proof Lets consider $d(e_{ij})= \displaystyle \sum_{x \in L_i}C^{ii}_{xi}e_{xj} + C_{ij}^{ij}e_{ij} + \sum_{y \in R_j}C^{jj}_{jy}e_{iy}$
for all $e_{ij} \in B$, where the coefficients $C^{ij}_{xy}$ are subject to the following relations
$$C^{jj}_{jk} + C^{kk}_{jk} = 0,  \ \ \ \ if j \leq k;$$
$$C^{ij}_{ij} + C_{jk}^{jk} = C_{ik}^{ik} , \ \ \ \ \ if i \leq j, j \leq k.$$ By \cite[Theorem 2.2]{Xiao} $d$ is a derivation.
First we check that
\begin{eqnarray}\label{id7}
\phi(e_{ij}e_{kl}) = \phi(e_{ij})e_{kl} + e_{ij}d(e_{kl}),
\end{eqnarray}
for all $e_{ij}, e_{kl} \in B$. We split the argument into two cases.

\noindent Case 1: $j \neq k$. Since $\phi(e_{ij}e_{kl}) = 0$, it suffices to prove that $\phi(e_{ij})e_{kl} + e_{ij}d(e_{kl}) = 0$.
By (\ref{id6}) we get 
\begin{eqnarray*} 
\phi(e_{ij})e_{kl} + e_{ij}d(e_{kl}) &=& (\Xi(e_{ii})e_{ij} + e_{ii}\tau(e_{ij}))e_{kl} + e_{ij}d(e_{kl}) \\&=& e_{ii}\tau(e_{ij})e_{kl} + e_{ij}d(e_{kl}).
\end{eqnarray*}
If $i \neq k$ then
\begin{eqnarray*}
e_{ii}\tau(e_{ij})e_{kl} + e_{ij}d(e_{kl}) &=& e_{ii}\tau(e_{ij})e_{kl} + e_{ij}d(e_{kk})e_{kl} \\&=& e_{ii}(\tau(e_{ij})e_{kk} + e_{ij}d(e_{kk}))e_{kl} \\&=& e_{ii} 0 e_{kl} \\ 
&=& 0,
\end{eqnarray*}
by Lemma \ref{xiao1} and $\tau(e_{ij})e_{kk} = \tau(e_{ij}e_{kk}) - e_{ij}\tau(e_{kk})$.
Finally, if $i =k$, then
\begin{eqnarray*}
e_{ii}\tau(e_{ij})e_{il} + e_{ij}d(e_{il}) &=& e_{ii}\tau(e_{ij})e_{il} + e_{ij}d(e_{ii}e_{il}) \\&=& e_{ii}\tau(e_{ij})e_{il} + e_{ij}d(e_{ii})e_{il} \\&=& (e_{ii}\tau(e_{ij}) + e_{ij}d(e_{ii})e_{il} \\&=& (\tau(e_{ij}) - \tau(e_{ii})e_{ij} \\&-& \tau(e_{ij})e_{ii} - e_{ij}\tau(e_{ii}) + e_{ij}d(e_{ii}))e_{il} \\&=& e_{ij}(d(e_{ii}) - \tau(e_{ii}))e_{il} = 0. 
\end{eqnarray*}

\noindent Case 2: $j = k$. We must prove that
$$\phi(e_{il}) = \phi(e_{ij})e_{jl} + e_{ij}d(e_{jl}).$$
Assume $i < j < l$. As a consequence of (\ref{id6}), 
\begin{eqnarray*}
\phi(e_{ij})e_{jl} + e_{ij}d(e_{jl}) &=& (\Xi(e_{ii})e_{ij} + e_{ii}\tau(e_{ij}))e_{jl} + e_{ij}d(e_{jl}) \\&=& \phi(e_{il}) - e_{ii}(\tau(e_{il}) - \tau(e_{ij})e_{jl} - e_{ij}d(e_{jl}))\\&=& \phi(e_{il}) - e_{ii}(e_{ij}\tau(e_{jl}) + \tau(e_{jl})e_{ij} + e_{jl}\tau(e_{ij}) \\&-& e_{ij}d(e_{jl})) \\&=& \phi(e_{il}) - e_{ij}(\tau(e_{jl})- d(e_{jl})) = \phi(e_{il}). 
\end{eqnarray*}
If $i=j <l$, then 
\begin{eqnarray*}
\phi(e_{ii})e_{il} + e_{ii}d(e_{il}) &=& \Xi(e_{ii})e_{il} + e_{ii}\tau(e_{il}) + e_{ii}d(e_{il}) - e_{ii}\tau(e_{il}) \\&=& \Xi(e_{ii})e_{il} + e_{ii}\tau(e_{il}) = \phi(e_{il}).
\end{eqnarray*}
If $i < j = l$ then 
\begin{eqnarray*}
\phi(e_{ij})e_{jj} + e_{ij}d(e_{jj}) &=& (\Xi(e_{ii})e_{ij} + e_{ii}\tau(e_{ij}))e_{jj} + e_{ij}d(e_{jj}) \\&=& \Xi(e_{ii})e_{ij} + e_{ii}\tau(e_{ij})+ e_{ii}\tau(e_{ij})e_{jj} \\&-& e_{ii}\tau(e_{ij}) + e_{ij}d(e_{jj}).
\end{eqnarray*}
Since $e_{ii}\tau(e_{ij})e_{jj} = C_{ij}^{ij}e_{ij}$, $e_{ii}\tau(e_{ij}) = C_{ij}^{ij}e_{ij} + \displaystyle \sum_{y \in R_j}C^{jj}_{jy}e_{iy}$ and $e_{ij}d(e_{jj}) = C_{jj}^{jj}e_{ij} + \displaystyle \sum_{y \in R_j}C^{jj}_{jy}e_{iy}$ it follows that $e_{ii}\tau(e_{ij})e_{jj} - e_{ii}\tau(e_{ij}) + e_{ij}d(e_{jj}) = 0$.
Hence $\phi(e_{ij})e_{jj} + e_{ij}d(e_{jj}) = \Xi(e_{ii})e_{ij} + e_{ii}\tau(e_{ij}) = \phi(e_{ij})$.
If $i=j=l$, by (\ref{id3}) we obtain $\phi(e_{ii}) = \Xi(e_{ii}) = \Xi(e_{ii})e_{ii} + e_{ii}\tau(e_{ii}) = \phi(e_{ii})e_{ii} + e_{ii}d(e_{ii})$.
Thus, for all $e_{ij}, e_{kl} \in B$, we get $\phi(e_{ij}e_{kl}) = \phi(e_{ij})e_{kl} + e_{ij}d(e_{kl})$. 
Finally, linearity of $\phi$ yields $\phi(ab) = \phi(a)b + ad(b)$ for all $a,b \in I(X,\R)$, which proves that $\phi$ is a generalized derivation. \ep
\\ 

We are now in a position to prove the main result of this paper.

\begin{thm}\label{main}
Let $\R$ be a $2$-torsion free commutative ring with identity. Then any generalized Jordan derivation of the incidence algebra $I(X,\R)$ is
a generalized derivation.\end{thm}
\proof Put $\Psi = \Xi - \phi$, then $\Psi(e_{ij}) = \Xi(e_{ij}) - \phi(e_{ij})$ and $\Psi(e_{ii}) = \Xi(e_{ii}) - \phi(e_{ii}) = 0$ for all $e_{ii} \in B$.
Since $\Psi$ is a generalized Jordan derivation then $\Psi(e_{ij}) = \Psi(e_{ij}e_{jj} + e_{jj}e_{ij}) = \Psi(e_{ij})e_{jj} + \Psi(e_{jj})e_{ij} = \Psi(e_{ij})e_{jj}$.
According to  (\ref{id4}) and (\ref{id6}), if $i < j$ we have
\begin{eqnarray*}
\Psi(e_{ij}) &=& \Xi(e_{ij})e_{ii} + e_{ij}\tau(e_{ii}) \\&=& (\phi(e_{ij}) + \Psi(e_{ij}))e_{ii} + e_{ij}\tau(e_{ii}) \\&=& \phi(e_{ij})e_{ii} + e_{ij}\tau(e_{ii}) +\Psi(e_{ij})e_{ii} \\&=& \phi(e_{ij}e_{ii}) + \Psi(e_{ij})e_{ii} \\&=& \Psi(e_{ij})e_{ii}. 
\end{eqnarray*}
Thus $\Psi(e_{ij}) = \Psi(e_{ij})e_{jj} = 0$.
Therefore $\Psi = \Xi - \phi = 0$ and $\Xi$ is a generalized derivation. \ep

\vspace{0.5cm}

As a consequence of our Theorem we have the following result.

\begin{cor}[Theorem 3.3 \cite{Xiao}] 
Let $\R$ be a $2$-torsion free commutative ring with identity. Then every Jordan derivation of the incidence algebra $I(X,\R)$ is
a derivation.
\end{cor}

\end{document}